\documentclass[11pt,psamsfonts]{amsart}

\usepackage{amssymb,amscd,amsmath,amsthm}

\chardef\bslash=`\\ % p. 424, TeXbook
\newtheorem{theorem}{Theorem}[section]
\newtheorem{corollary}[theorem]{Corollary}
\newtheorem{lemma}[theorem]{Lemma}
\newtheorem{proposition}[theorem]{Proposition}

\theoremstyle{remark}
\newtheorem{remark}[theorem]{Remark}

\theoremstyle{definition}

\numberwithin{equation}{section}

\newcommand{\thmref}[1]{Theorem~\ref{#1}}

\newcommand{\proref}[1]{Proposition~\ref{#1}}
\newcommand{\lemref}[1]{Lemma~\ref{#1}}
\newcommand{\corref}[1]{Corollary~\ref{#1}}

\newcommand{\clsp}{\overline{\operatorname{span}}}
\newcommand{\Ind}{\operatorname{Ind}}
\newcommand{\sign}{\operatorname{sign}}
\newcommand{\nat}{\mathbb N^{\times}}
\newcommand{\inv}{^{-1}}
\newcommand{\N}{\mathbb N}
\newcommand{\Z}{\mathbb Z}
\newcommand{\Q}{\mathbb Q} 
\newcommand{\C}{\mathbb C}
\newcommand{\R}{\mathbb R}
\newcommand{\PP}{\mathbb P}
\newcommand{\primes}{\mathcal P}
\newcommand{\adeles}{\mathcal A}
\newcommand{\bcheck}{\mathcal C_{\Q}}
\newcommand{\qmz}{\Q / \Z}
\newcommand{\qmzhat}{\widehat{\qmz}}
\newcommand{\qat}{\Q^*_+}
\newcommand{\qhat}{\widehat{\qat}}

\newcommand{\Prim}{\operatorname{Prim}}

\title[The Hecke $C^*$-algebra of Bost and Connes]{The ideal structure of the Hecke
$C^*$-algebra\\ of Bost and Connes}
\date{October 11, 1999.}
\thanks{This research was supported by the Australian Research
Council.} 

\author[Marcelo Laca]{Marcelo Laca}
\author[Iain Raeburn]{Iain Raeburn}
\address{Department of Mathematics \\
    The University of Newcastle\\  NSW  2308\\ AUSTRALIA}
\email{marcelo { and} iain @math.newcastle.edu.au}

\begin{document}
\begin{abstract} We compute explicitly the primitive ideal space of the
Bost-Connes Hecke $C^*$-algebra by embedding it as a full corner in a  transformation group
$C^*$-algebra and applying a general theorem of Williams.  This requires
the computation of the quasi-orbit space for the action of $\Q^*_+$  on
the space $\adeles_f$ of finite adeles. We then carry out a similar
computation for the action of $\Q^*$ on the space $\adeles = \adeles_f \times \R$ of full adeles.
\end{abstract}
\maketitle
\section*{Introduction}
 
 In their work on phase transitions in number theory, Bost and Connes
introduced a  noncommutative Hecke $C^*$-algebra $\bcheck$ and gave a
presentation of this algebra. The present authors recognised this
presentation as that of a semigroup crossed product 
$C^*(\qmz) \rtimes_\alpha \N^\times$, and analysed the representations
of $\bcheck$ using techniques previously developed for studying
Toeplitz algebras \cite{bcalg}. Laca showed in \cite{diri} that this
approach gives useful insight to many of the constructions of Bost and
Connes, and used the universal properties of semigroup crossed products
to simplify and extend their arguments.

From a $C^*$-algebraic point of view, however, the Hecke $C^*$-algebra
$\bcheck$ remains poorly  understood. Here we shall contribute a
description of the ideal structure of $\bcheck$. More precisely, we
shall describe the primitive ideal space $\Prim \bcheck$ and the
topology on it: since the lattice of ideals in any $C^*$-algebra $A$  is
isomorphic to the lattice of closed subsets of
$\Prim A$, this gives a complete description of the ideal structure.

To do this, we use the realisation of $\bcheck$ as 
$C^*(\qmz) \rtimes_\alpha \N^\times$. Following a standard route
(see \cite{mur,minautex}), we
dilate the endomorphisms $\alpha_n$ to automorphisms $\alpha_n^\infty$
of a direct limit $C^*(\qmz)_\infty$, and then realise 
$C^*(\qmz) \rtimes_\alpha \N^\times$ as a full corner in the ordinary
crossed product 
$C^*(\qmz)_\infty \rtimes_{\alpha^\infty} \Q^*_+$. The direct limit
$C^*(\qmz)_\infty$ is naturally isomorphic to the algebra
$C_0(\adeles_f)$ of continuous functions on the space $\adeles_f$ of
finite adeles, and the action
$\alpha^\infty$ is conjugate to the action of the diagonally embedded
copy of $\Q^*_+$ by division on $\adeles_f$. We analyse
$\Prim(C_0(\adeles_f) \rtimes \Q^*_+)$ using a theorem of Williams which
identifies it as a quotient space of $\mathcal Q \times \qat$, where
$\mathcal Q$ is the quasi-orbit space for the action (see \S1). The main
part of our work, therefore, is to compute this quasi-orbit space; the
parametrisation we obtain allows us to give an elegant description of
the topology (see \thmref{summary}).

The space $\adeles_f$ of finite adeles
tells only part of the story, and one should also consider
the ``infinite place" --- that is, use the space $\adeles = \adeles_f
\times \R$ of full adeles. (In Connes' latest analysis \cite{con-rzf},
for example, he has implicitly discarded $\bcheck$ in favour of  $C_0(\adeles)\rtimes\Q^*$.)
 We have been able to extend our analysis of the quasi-orbits to the action of $\Q^*$ on 
$\adeles$, and thus obtain a
parametrisation of $\Prim(C_0(\adeles) \rtimes \Q^*)$ and a
description of the topology. The inclusion of the
infinite place significantly changes the  primitive ideal space,
largely because the orbits of invertible adeles are now closed rather than
dense, giving rise to a copy of the idele class group inside 
the primitive ideal space (see \corref{Primfulladeles}).

\section{Williams' theorem} Consider an action of a second
countable locally compact abelian group $G$ on a  second countable
locally compact space $X$. Recall that a  quasi-orbit for the action is
the set of points with a given orbit closure; the {\em quasi-orbit
space} $\mathcal Q(X / G)$
 is the quotient of $X$ by the relation 
$$ x \sim y \quad \iff \quad \overline{G\cdot x} = \overline{G\cdot
y}.
$$ Notice that, because $G$ is abelian, the isotropy groups
$G_x$ are constant on quasi-orbits. The Mackey machine, as presented by
Green \cite{gre}, says that every primitive ideal is associated to a
unique quasi-orbit (in a sense made precise in \cite{gre}) and is
induced from the isotropy group of that quasi-orbit. To get specific
representations, let $\varepsilon_x $ denote the representation 
$f \mapsto f(x)$ of $C_0(X)$. Then for every $\gamma \in \hat G$, 
$(\varepsilon_x,\gamma|_{G_x})$ is a covariant representation of
$(C_0(X), G_x)$, and the induced representation
$\Ind_{G_x}^G(\varepsilon_x\times \gamma)$ is an irreducible
representation of $C_0(X) \rtimes G$ (see \cite[Proposition 4.2]{wil}). We shall use the following
form of Williams' theorem.

\begin{theorem}[Williams]\label{williamsthm}  Define an equivalence
relation on $\mathcal Q(X / G)\times \hat{G}$ by
$$ ([x],\gamma) \sim ([y],\chi)  \quad \iff \quad [x] = [y]
\text{\  and \ } \gamma|_{G_x} = \chi|_{G_x},
$$ where $G_x$ denotes the isotropy group at $x$.
 Then the map 
\begin{equation}\label{quasiorbit}
 ([x], \gamma) \mapsto \ker \Ind_{G_x}^G(\varepsilon_x\times \gamma)
\end{equation} induces a homeomorphism of the quotient 
$(\mathcal Q(X / G)\times \hat{G})/_\sim $ onto $\Prim (C_0(X) \rtimes
G)$.
\end{theorem}

\begin{proof}  An examination of the proof of \cite[Lemma
4.10]{wil} shows that the map $(x,\gamma) \mapsto \ker \Ind
(\varepsilon_x\times \gamma)$ of $X\times \hat{G}$ into $\Prim
(C_0(X)\rtimes G)$ is the map discussed at the beginning of
\cite[\S5]{wil}. Since second countable transformation groups are
automatically quasi-regular (\cite[\S5]{eff-hah}; see \cite[Corollary
19]{gre} for  a more general result), it follows from
\cite[Corollary 3.2]{goo-ros} that every second countable
transformation group $(X,G)$ with $G$ amenable is Effros-Hahn regular
in the sense of 
\cite[Definition 4.12]{wil}. This applies in particular when $G$ is
abelian, so the  present theorem is a restatement of \cite[Theorem
5.3]{wil} as it applies to second countable abelian transformation groups; that
Williams' quotient $\Gamma$ of $X\times \hat G$ can be viewed as a
quotient of $\mathcal Q(X/G) \times \hat G$ is observed before
\cite[Theorem 5.3]{wil}.
\end{proof}

\begin{remark}\label{Macrosslambda} In our case, the isotropy groups
will be either $G$ or $\{e\}$. When $G_x = G$, the induced
representation $\Ind_{G_x}^G(\varepsilon_x\times \gamma)$ is just
$\varepsilon_x\times \gamma$. When $G_x = \{e\}$, the representation 
$\Ind_{G_x}^{\{e\}}\varepsilon_x$ is equivalent to $M_x \times
\lambda$ acting on $L^2(G)$, where
$$ (M_x(f)\xi)(s) = f(s\cdot x)\xi(s) \qquad \text{and} \qquad 
\lambda_t (\xi)(s) = \xi(t\inv s).
$$ To see this, note that the representation considered in \cite{wil}
is given by the left action $(R,V)$ of
$(C_0(X), G)$ on the Hilbert $C_0(X)$-module $\overline{C_c(G,X)}$
described  in the middle of page 340 of \cite{wil}. The unitary $U:
\overline{C_c(G,X)} =\overline{C_c(G,X)} \otimes \C \to L^2(G)$ of
\cite[Lemma 2.14]{wil} is given by $U_z(r) = z(r,r\cdot x)$, and the
calculations
$$ U(R_x(f) z)(r) = (R_x(f) z) (r,r\cdot x) = f(r\cdot x) z(r,r\cdot x)
= (M_x (f) U(z) ) (r)
$$ and
$$ U(V(t)z)(r) = (V(t)z)(r,r\cdot x) = z(t\inv r,t\inv r \cdot x) =
(\lambda_t U(z))(r)
$$ show that $U$ intertwines $\Ind_{G_x}^{\{e\}}\varepsilon_x$
with $M_x \times
\lambda$.
\end{remark}

\section{The primitive ideal space of $\bcheck$}\label{prim-bcalg}

Let $\primes$ be the set of   prime numbers and, for $p\in \primes$, 
let  $\Q_p$ be the  field of 
$p$-adic numbers, obtained by completing the rationals under the
$p$-adic absolute value. The ring of $p$-adic integers $\Z_p$  is the
closure of the integers in $\Q_p$, and is a compact open subring of
$\Q_p$. The ring $\adeles_f$ of finite adeles is by definition  the
restricted product 
$$
\prod_{p\in \primes} (\Q_p; \Z_p) = \Big\{ (a_p) \in \prod_{p\in
\primes} \Q_p:
 a_p \in \Z_p \text{ for all but finitely many } p \in \primes\Big\},
$$  on which the topology is defined by stating that 
 a neighbourhood base at the identity is given by the collection of sets
$\prod_F V_p \times \prod_{p\notin F}\Z_p $, with $F$ a finite subset
of $\primes$ and
 $V_p$  a neighbourhood of $ 0$ in $ \Q_p$.
There is a canonical diagonal embedding of $\Q$ in $\adeles_f$, defined
via the canonical embeddings of $\Q$ in $\Q_p$ for each
$p$. Notice that $r \in \Q$ embeds as a $p$-adic integer for all $p$
relatively prime to the reduced denominator of $r$. The group $\Q^*_+$
of positive rationals embeds as units in $\adeles_f$. For these
basic number theoretic facts we refer to \cite{ram-val,tate}.

By \cite[Corollary 2.10]{bcalg}, the Bost-Connes $C^*$-algebra $\bcheck$
is isomorphic to a semigroup crossed product 
$C^*(\qmz)\rtimes_\alpha \nat$, which by \cite[Proposition
3.2.1]{minautex}   is isomorphic to a full corner in the transformation
group $C^*$-algebra
 $C_0(\adeles_f)\rtimes \qat$ for the
 action of the  positive rationals $\qat$ by division on $\adeles_f$.
Since a full corner has the same representation theory as its
ambient algebra (see \lemref{morita-eq} below for the specific
statement), it will suffice to compute the primitive ideal space of $
C_0(\adeles_f)\rtimes
\qat$, which we can do using  Williams' theorem.  We need to compute
the orbit closures and the duals of the isotropy groups at each point.
The isotropy groups are easy to compute:

\begin{lemma}\label{isotrop} The isotropy group at every $a \in
\adeles_f \setminus \{0\}$ is trivial,
 and the isotropy group at $0 \in \adeles_f$ is  $\qat$.
\end{lemma}
\begin{proof} For each $r\in \qat$, we have $r  a := (ra_p)_{p\in
\primes}$, so 
$r   a = a $ implies that $r a_p = a_p$ in $\Q_p$  for every prime $p$.
Clearly $a \neq 0$ implies that $a_p \neq 0 $ for at least one prime
$p$, and, for such a prime, $r a_p = a_p $ implies $r = 1$, because
$\Q_p$ is a field. The second assertion is obvious.
\end{proof}

\begin{remark}\label{quot>dirsum} Because the isotropy group at $0$ is all of $\qat$,
Williams' theorem says that there must be a copy of $\qhat$ in the 
primitive ideal space of $\bcheck$.  It is easy to identify this copy:
 if $\epsilon :  C^*(\Q/\Z) \to \C$ is the augmentation  homomorphism
given by $\epsilon (\delta_r) = 1$ for every $r\in \qmz$, then  for
each $\gamma \in \qhat$,
 the pair $( \epsilon, \gamma|_{\nat})$ is a one-dimensional covariant
representation of  $(C^*(\Q/\Z), \nat, \alpha)$, giving a
one-dimensional representation 
 of $C^*(\Q/\Z) \rtimes_\alpha \nat$.
\end{remark}
 Next we need to compute the closure of an orbit. To
simplify things we note that the orbits for the action of $\qat$ by
division are  the same as for the action by multiplication, and compute
the latter. 

\begin{lemma}\label{binqa}
 For each $a \in \adeles_f$, the closure of the corresponding orbit is
\begin{equation}\label{binqaeq}
\overline{\qat   a} = \{ b \in \adeles_f: a_p = 0 \Longrightarrow b_p
=0\}.
\end{equation}
\end{lemma}

\begin{proof} It is clear that $\subset$ holds, so we concentrate on
proving $\supset$. Let
$b$ belong to the right-hand side, and notice that it is enough to show
that
$mb \in \overline{\qat   a}$ for some integer $m$. We may choose $m$
such that $mb \in \prod_\primes \Z_p$, and then we can easily write
down a typical basic open neighbourhood of $mb$: 
$$ V = \prod_{p\in F} V_p \times \prod_{p\notin F} \Z_p
$$  for $F$ a finite set of primes and  $V_p$ open sets in $\Z_p$.  We
can further assume that $a \in \prod \Z_p$ without changing its orbit
or the set of primes $p$  for which $a_p$ vanishes.

 Let $F' := \{p \in F: a_p \neq 0\}$, and notice that for every $p\in
F'$ we can write $a_p = p^{v_p(a_p)}u_p$ for some unit $u_p \in
\Z_p^*$. Then $ p^{v_p(a_p)}a_p\inv V_p = u_p\inv V_p$ is an open set
in $\Z_p$, and hence contains a positive integer $k_p$; we then have
that
$k_p p^{-v_p(a_p)}\in a_p\inv V_p$.

We now choose $l_p$ large enough to ensure that the ball $B(k_p
p^{-v_p(a_p)}; p^{-l_p + v_p(a_p)})$ in $\Q_p$ is contained in $a_p
\inv V_p$. (Notice that because $V_p$ is open in
$\Z_p$, and $a_p$ is nonzero in $\Q_p$, we have that $a_p\inv V_p$ is
open in $\Q_p$.) By the Chinese Remainder Theorem, there is an integer
$n$ such that
$$
 n \equiv k_p \Big(\prod_{q\in F' \setminus \{p\}} q^{v_q(a_q)}\Big)
\mod p^{l_p}
\qquad \text{ for all } p\in F'.
$$
 But now for each $p\in F'$,
$$
\left| \frac{k_p}{p^{v_p(a_p)}} - \frac{n}{\prod_{q\in F'}
q^{v_q(a_q)}}\right|_p = 
\left| \frac{(k_p \prod_{q\in F'\setminus\{p\}} q^{v_q(a_q)}) -
n}{\prod_{q\in F'} q^{v_q(a_q)}}\right|_p
\leq p^{-l_p + v_p(a_p)}.
$$  Next define 
$$  r := \Big(\frac{n}{\prod_{q\in F'} q^{v_q(a_q)}}\Big),
$$ and notice that  for all $p\in F'$ we have
$r  \in  B(k_p p^{-v_p(a_p)}; p^{-l_p + v_p(a_p)}) \subset a_p\inv
V_p$, so that
$r a_p\in V_p$. If $p\in  F\setminus F'$, in which case $a_p = 0 $,
then $ra_p = 0 = mb_p$, so clearly $ra_p \in  V_p$. If $p \notin F$,
then  $ {\prod_{q\in F'} q^{v_q(a_q)}} $ is coprime to $p$, and $r$
embeds as a $p$-adic integer (that is, $r \in \Z_p$), and hence $ra_p
\in \Z_p$. Thus 
$r a \in  \prod_{p\in F}V_p \times
\prod_{p\in \primes \setminus F} \Z_p$, and $mb \in
\overline{\qat   a}$.
\end{proof}

The quasi-orbit space $\mathcal Q(\adeles_f/\qat)$ is the quotient
  of $\adeles_f$ by the equivalence
relation $\sim$ defined by
$ a\sim b $ if and only if $ \overline{\qat   a} = \overline{\qat b}$.
The above characterisation of the orbit closures shows that $a\sim b$
if and only if $\{p: a_p = 0 \} = \{p: b_p = 0 \}$, hence
 the quasi-orbits are in one-to-one correspondence with the subsets of
the set of prime numbers. The power set $2^\primes$ has a natural
{\em power-cofinite topology} in which the basic open sets are
$$  U_G := \{ T\in 2^\primes :  T \cap G = \emptyset\}
 =
 2^{ \primes \setminus G},
$$
 where $G$ ranges over the finite subsets of $\primes$. 
 This collection is indeed the basis for a topology because $U_G \cap
U_H = U_{G\cup H}$.

\begin{proposition}\label{qospace} For each $a\in \adeles_f$ define
$ S(a) := \{ p\in \primes: a_p = 0\}$. Then the map $a \mapsto S(a)$
induces a homeomorphism of the quasi-orbit space 
$\mathcal Q(\adeles_f/\qat)$ onto the space $2^\primes$ with the
power-cofinite topology defined above.
\end{proposition}

\begin{proof} Let $q: a \mapsto q(a)$ be the quasi-orbit map of
$\adeles_f $ to $ \mathcal Q(\adeles_f/\qat)$. It follows easily from
the characterisation of orbit closures that the map $a \mapsto S(a)$
factors through  $q$, and that $q(a) \mapsto S(a)$ is a bijection.

By  the Lemma on page 221 of \cite{gre} we know that the map
$q:\adeles_f \to \mathcal Q(\adeles_f/\qat)$ is continuous and open.
Hence the sets $q(V)$, as $V$ runs through a basis for the topology on $\adeles_f$, form a basis for
the quotient topology on $\mathcal Q(\adeles_f/\qat)$. By definition, a typical basic
open neighbourhood of $a$ in $\adeles_f = \prod_{p\in \primes}
(\Q_p:\Z_p)$ is  $W = a + V$ with $V$ a product of the form
$$  V = \prod_{p\in F} V_p \times \prod_{p\notin F} \Z_p,
$$ in which $F$ is a finite subset of $\primes$ and $V_p$ is an open
neighbourhood of $0$ in $\Q_p$ for every $p\in F$. After relabelling $a
+ \Z_p$ as $V_p$ for those (finitely many) $p \notin F$ for which
$a_p \notin \Z_p$, and enlarging the set $F$ accordingly, we see that
$W$ is a product of the same form as $V$. 

Since the bijection $q(a) \mapsto S(a)$ carries
$q(V)$ to $S(V)$, we want to prove that  the collection $\{S(V)\}$ with
$V$ as above is precisely the basis $\{U_G\}$ for the topology on $2^\primes$. Given
$V$ as above, the set $ \{ p\in F:  0 \notin V_p\}$ is finite, and every finite subset of $\primes $
arises for some $V$. Thus it suffices to prove that for every $V$ as above,
$$  S(V) =  U_{\{ p\in F:  0 \notin V_p\}}.
$$
 First we prove $\subset$:
\begin{eqnarray*} b \in V  & \Longrightarrow  & b_p \in V_p \quad
\text{ for all } p\in F \\ 
 &\Longrightarrow & b_p \neq 0  \quad \text{ if } 0 \notin V_p \\
 &\Longrightarrow &  S(b)  \cap \{ p\in F:  0 \notin V_p\} = \emptyset\\
 &\Longrightarrow & S(b) \in U_{\{ p\in F:  0 \notin V_p\}}.
\end{eqnarray*}
To prove $\supset$, let $G : =\{ p\in F:  0 \notin V_p\}$ and suppose
$T$ is a subset of $\primes$ such that $T \cap G = \emptyset$.  We need
to find $b = (b_p) \in V$ such that $S(b) = T$. 
If $p\in T$ choose $b_p = 0$;  if $p \notin T$ but $p\in F$, choose
$b_p \in V_p \setminus\{0\}$; finally, if $p \notin T \cup F$, simply
take $b_p = 1 \in \Z_p$. It is clear that $S(b) =T$. To see that $b \in
V$, notice that $b_p \in V_p$ for $p \in F \cap T$ because $T\cap G
=\emptyset$, and obviously $b_p \in V_p$ for $p \in F \setminus T$.
\end{proof}

Now that we have described the quasiorbit space and the isotropy groups, our description of $\Prim
(C_0(\adeles_f)
\rtimes
\qat)$ follows immediately from
\thmref{williamsthm}:

\begin{proposition} \label{combine} Define a relation on $2^\primes
\times \qhat$ by
$$ (S,\gamma) \sim (T,\chi) \quad \iff \quad \left\{
\begin{array}{ll} 
                S = T & \text{ if } S \neq \primes \\
                S= T \text{ and } \gamma = \chi  & \text{ if } S =
\primes.
              \end{array} \right.
$$ Then $\sim$ is an equivalence relation. For each nonzero $a\in
\adeles_f $ the ideal $I_S:= \ker
\Ind \varepsilon_a = \ker(M_a \times \lambda)$ depends only on $S :=
\{p: a_p = 0\}$, and 
 the maps
$$ S \mapsto I_S \quad \text{ and }\quad \gamma \mapsto
\ker(\varepsilon_0\times \gamma)
$$ combine to give a homeomorphism of 
$(2^\primes  \times \qhat) /_\sim $ onto $\Prim (C_0(\adeles_f) \rtimes
\qat)$.
\end{proposition}
\begin{proof} This follows from a direct application of
\thmref{williamsthm},  using the characterisation of the quasi-orbit
space given by \proref{qospace}.
\end{proof}

\begin{remark} \label{disjunion}
Since the isotropy groups are either trivial or $\qat$,
we can also describe 
$(2^\primes  \times \qhat) /_\sim $ as the disjoint union
$$
 (2^{\primes} \setminus \{\primes\}) \sqcup {\qhat}
$$
 topologised as follows: the subspace
 $2^{\primes} \setminus \{\primes\}$ is open and has the power-cofinite topology, and the subspace
$\qhat$ has its usual compact topology. Moreover, all of $2^{\primes} \setminus
\{\primes\}$  is contained in every neighbourhood of any point in
$\qhat$.
\end{remark}

So far we have been working with the dilated system consisting of the
action of $\qat$ on $\adeles_f$ rather than with $\bcheck$. The algebra
$\bcheck$ is isomorphic to the corner in $C_0(\adeles_f) \rtimes \qat$
corresponding to the projection $i_{C_0(\adeles_f)}(1_\mathcal Z)$,
where $\mathcal Z$ is the maximal compact subring of $\adeles_f$ given
by $\mathcal Z =
\prod_{p \in \primes} \Z_p$ (see \cite[\S3]{minautex}).  To pass to
representations of
$\bcheck$ we shall use the following standard lemma.

\begin{lemma}\label{morita-eq} Let $B$ be a $C^*$-algebra and $e\in M(B)$
a full projection. Then the map which sends a representation $\pi $ of
$B$ to its compression $\pi|_{eBe}$ to $\pi(e) \mathcal H_\pi$ respects
unitary equivalence and induces a homeomorphism of $\Prim B$ onto
$\Prim eBe$.
\end{lemma}

\begin{proof} The space $eB$ is an $eBe$--$B$ imprimitivity bimodule.
If $\pi$ is a representation of $B$, the corresponding representation
$eB\text{--}\Ind_B^{eBe} \pi$ acts on $eB \otimes_B
\mathcal H_\pi$ via 
$$ (eB\text{--}\Ind_B^{eBe} \pi) (e c e) (eb \otimes h) = e c e b
\otimes h.
$$ But the map 
$ U:eB \otimes_B \mathcal H_\pi \to \pi(e)\mathcal H_\pi$ defined by
$U(eb\otimes h) =
\pi(eb)h$ is unitary, and intertwines the natural left action of $eBe$
with the restriction of $\pi$. Thus $eB\text{--}\Ind_B^{eBe} \pi$ is
equivalent to the compression of $\pi|_{eBe}$ to its essential subspace
$\pi(e) \mathcal H_\pi$. The result therefore follows from standard
properties of Morita equivalence \cite[Corollary 3.33]{rae-wil}.
\end{proof}

Let $S$ be a proper subset of $\primes$,  and choose $a\in \adeles_f$
such that 
$a_p = 0$ if and only if $p\in S$. Then, using the notation of
Remark~\ref{Macrosslambda},
$$
\ker \Ind^{\qat}_{\{e\}} \varepsilon_a =  \ker (M_a \times \lambda)
$$
is the primitive ideal of  $C_0(\adeles_f)\rtimes \qat$ corresponding
to $S$. 
 The corresponding representation of $\bcheck$ is the compression of 
$(M_a \times \lambda)|_{\bcheck}$ to the subspace 
$$ M_a(1_{\mathcal Z})\ell^2(\qat) =
 \clsp \{e_r : r\in \qat, \  r  a \in \mathcal Z\}
$$ and is determined by the covariant representation $(\pi_a, V)$ of
$(C^*(\qmz), \nat, \alpha)$
 where $V_n e_r = e_{nr}$ and $\pi_a(f) e_r = \hat{f}(r  a) e_r$. Here
the Fourier transform $\hat f$ of $f \in C^*(\qmz)$ is viewed as a
function in $C(\mathcal Z)$ using the canonical identification of $\mathcal Z$ with 
$\qmzhat$ (see, for example, \cite[page 360]{diri}).  We summarise the results in the following
theorem.

\begin{theorem}\label{summary} For each proper subset $S$ of $\primes$,
choose $a\in \mathcal Z= \prod_p \Z_p$ such that $S=\{p\in \primes : a_p= 0\}$.
 Then the representation $\pi_a\times V$ is irreducible and its kernel
depends only on $S$.
 The maps
$$ S \mapsto \ker (\pi_a \times V) \quad \text{ and } \quad 
\gamma \mapsto \ker(\varepsilon_0 \times \gamma)
$$  combine to give a homeomorphism  of $  (2^\primes \times
\qhat)/_\sim$ onto
$\Prim \bcheck$. 
\end{theorem}

\begin{proof} The preceding discussion shows that the
map of $(2^\primes  \times \qhat) /_\sim $ onto $\Prim \bcheck$ is the
composition of the homeomorphism of $(2^\primes  \times \qhat) /_\sim $
onto $\Prim (C_0(\adeles_f) \rtimes \qat)$  of \proref{combine} with
the homeomorphism of
\lemref{morita-eq}; in particular, it is a homeomorphism.
\end{proof}

 This shows that for $a \neq 0$, the representations $\pi_a\times V$
and $\pi_b \times V$ are weakly equivalent (have the same kernel)
if and only if $a$ and $b$ are in the same quasi-orbit. It is natural
to ask when the representations $\pi_a\times V$ and $\pi_b \times V$
are unitarily equivalent. First we need to characterise the zero
divisors of $\mathcal Z$.

\begin{lemma}\label{zerodiv}
 Let $a \in \mathcal Z =\prod_p \Z_p$.  Then the following are
equivalent:
\begin{enumerate}
\item[(1)] $a$ is not a zero divisor;
\item[(2)] $a_p \neq 0 $ for every $p \in \primes$;
\item[(3)] the set $\qat a \cap \mathcal Z$ is dense in $\mathcal Z$.
\end{enumerate}

\end{lemma}
\begin{proof} Since the operations in $\mathcal Z = \prod_p \Z_p$  are
componentwise, 
 (1) is trivially equivalent to (2).

Suppose (2) holds, let $F$ be a finite subset of $\primes$, and let 
$V_p$ be an open subset of $\Z_p$ for each $p\in F$, so that
$V = \prod_{p\in F} V_p \times \prod_{p \notin F} \Z_p$ is a basic open
set in
$\prod_p \Z_p$. Let $n = \prod_{p \in F}  p^{v_p(a_p)}$, and define
$ u = (u_p)_{p \in \primes}\in \prod_p \Z_p$ by
$$ u_p = \left\{ \begin{array}{ll} 
                \frac{1}{n} a_p & \text{ if } p\in F \\
                   1            & \text{ if } p\notin F.
              \end{array} \right.
$$  Then $u$ is a unit, $n$ divides $a$, and $\frac{1}{n}a_p = u_p$
for $p \in F$. Since multiplication by a unit is a homeomorphism,
$u^{-1} V$ is an open set, and since the canonical embedding of $\nat$
is dense in
${\mathcal Z}$, there exists $m\in \nat$ such that
$m \in u^{-1} V$. Thus $\frac{m}{n}a_p = m u_p \in V_p$ for every $p\in
F$, so $\frac{m}{n}a
\in V$, proving (3).

  Conversely, suppose $a_p = 0 $ for some $p$. Then $\qat a \cap
\mathcal Z$ is contained in the proper closed set $\{z: z_p = 0 \}$ and
cannot be dense, so (3) implies (2).
\end{proof}

\begin{proposition} Suppose $a\in \mathcal Z$ and $a \neq 0$.
\begin{enumerate}
\item[(1)] The representation $\pi_a\times V$ is faithful if and only
if $a $ is not a zero divisor. 

\item[(2)] Suppose $b\in \mathcal Z \setminus \{0\}$. The
representations
 $\pi_a\times V $ and $\pi_b \times V$ are unitarily equivalent if and
only if
 $a$ and $b$ are
in the same $\qat$-orbit.
\end{enumerate}
\end{proposition}
\begin{proof} By \cite[Proposition 3.7]{bcalg}, $\pi_a \times V$ is
faithful if and only if $\pi_a$ is faithful on $C^*(\qmz) \cong
C(\mathcal Z)$. Since the kernel of $\pi_a$ in $C(\mathcal Z)$ is
$\{g: g|_{\qat  a} = 0\}$, the assertion (1) follows from
\lemref{zerodiv}.

Next we prove (2). If $b = r a$ for some $r\in \qat$,  then $V_r$
implements the equivalence between $\pi_a\times V $ and $\pi_b \times
V$. Now suppose that $\pi_a\times V $ is unitarily equivalent to $\pi_b
\times V$; we have to show that $a$ and $b$ lie in the same
$\qat$-orbit.

Denote the canonical conditional expectation from $\bcheck$ to
$C^*(\qmz)$ by $\Phi$; then
$\omega_b: T \mapsto \widehat{\Phi(T)}(b)$ is the vector state $T
\mapsto \langle (\pi_b\times V) (T) e_1, e_1\rangle$. If $\pi_b \times
V $ is unitarily equivalent to $ \pi_a\times V$, then there is a unit
vector $\xi$ in the Hilbert space of $ \pi_a \times V$ such that
$\langle \pi_a\times V (T) \xi , \xi\rangle = \omega_b(T) $. Write $\xi
= \sum_{r}  c_r e_r$, where the sum is over the set $\{r \in
\qat: r a\in \mathcal Z\}$, and choose a finite subset  $F$ of $\qat$
such that $\sum_{r\notin F}  |c_r|^2 < 1/2 $. If $b \notin \qat a$
there exists $f\in C^*(\qmz)$ with $\|f\| \leq 1$ such that 
$\hat{f} (b) = 1$ and $\hat{f}(r a) =0$ for every $r\in F$. But then
$$ 1 = \hat{f}(b) = \omega_b(f)=  \langle \pi_a(f) \xi , \xi \rangle  
= \Big| \sum_r \langle \hat{f}(r a) c_r e_r,c_r e_r\rangle \Big|  \leq
\sum_{r\notin F}  |c_r|^2 < 1/2,
$$ which is a contradiction. So $b$ has to be in $\qat a$. 
\end{proof}

\section{The action of $\Q^*$ on the full adeles}\label{fulladeles}

In this section we compute the
primitive ideal space of the transformation group $C^*$-algebra $C_0(\adeles) \rtimes \Q^*$
of the multiplicative action of $\Q^*$ on the full adeles $\adeles = \adeles_f \times \R$.
If $a \in \adeles$, we write $a_f$ for the finite part $(a_p)_{p \in \primes}$ 
and $a_\infty$ for the value of $a$ at the infinite prime, so that $a = (a_f,a_\infty) \in
\adeles_f \times \R$. 

For each finite prime $p$, let $v_p: \Q_p \mapsto \Z \cup \{\infty\}$ be the $p$-adic valuation.
(If $x \in \Q_p^*$, then $v_p(x)$ is by definition the unique integer such that $x p^{-v_p(x)}$
is in $\Z_p^*$; if $x = 0$ we let $v_p(x) = \infty$.)
By abuse, if $a \in \adeles$, we write $v_p(a) = v_p(a_p)$. If $a$ is invertible, 
then $v_p(a) $ is finite for every $p$ and vanishes for almost all $p$; it follows that 
the product $\prod_{p\in \primes} p^{-v_p(a)}$ is a rational\ number, and it is easy to check that
$a \sign(a_\infty) \prod_{p\in \primes} p^{-v_p(a)} \in \prod_{p\in \primes} \Z_p^* \times \R^*_+$ 
for every invertible adele $a$.
Thus, every invertible adele has the form $ r u $ for some
$u \in \prod_p \Z_p^* \times \R^*_+$ and $r\in \Q^*$. Clearly every such product is invertible,
and the factorisation is unique, so we have a bijection
\begin{equation}\label{factor-ideles}
(r,u) \in \Q^* \times \big(\prod_{p\in \primes} \Z_p^* \times \R^*_+\big) \mapsto ru \in \adeles^*.
\end{equation}

As for the action of $\Q^*$ on $\adeles_f$, $0$ is the only adele with nontrivial 
isotropy, so the quasi-orbit space is the key to understanding the primitive 
ideal space. We begin by characterising the orbit closures.
\begin{lemma}\label{qorb-inv}
 If $u$ is an invertible adele, then the orbit $\Q^*  u$ is closed in $\adeles$.
\end{lemma}
\begin{proof} We
may assume that
$u = (u_f,u_\infty) \in \prod_p \Z_p^* \times \R^*_+$ without changing
the orbit. Suppose now $a $ is in the closure of $\Q^* u$, and let $r_n$
be a sequence of nonzero  rationals such that $r_n u \to a$ as $n \to
\infty$. Fix  $\epsilon < 1/(4u_\infty)$. Since $\prod_p \Z_p \times
(-\epsilon,\epsilon)$ is a neighbourhood of
$0$ in $\adeles$, there exists $N \in \N$ such that for every $n > N$,
$$ r_n u - a \in \prod_p \Z_p \times (-\epsilon,\epsilon).
$$ Then for every $m,n > N$ we have 
\begin{align*} (r_n  - r_m) u = (r_n u - a) - (r_m u - a)  &\in
\left(\prod_p \Z_p \times (-\epsilon,\epsilon)\right) - \left(\prod_p
\Z_p
\times (-\epsilon,\epsilon) \right)\\ & = \prod_p \Z_p \times
(-2\epsilon,2\epsilon).
\end{align*} 
Since $u_f \in \prod_p \Z_p^*$,   $(r_n - r_m) u
\in  \prod_p \Z_p $ implies
$r_n - r_m \in \Z$.  On the other hand, $(r_n - r_m)u_\infty \in
(-2\epsilon, 2\epsilon)$ implies
$$|r_n - r_m| < \frac{4\epsilon}{u_\infty} < 1.
$$  Since $r_n - r_m$ is an integer, the sequence
$r_n$ must be eventually constant, and hence $a = r_n u \in \Q^*
u$.
\end{proof}

Next we characterise the orbit closure of noninvertible adeles.
The result is similar to \lemref{binqa} for the finite
adeles, but the proof is slightly more involved because we need to 
control the size of the archimedean absolute value. 
  
\begin{lemma} \label{qorb-noninv}
If $a \in \adeles$ is not invertible, then 
$$
\overline{\Q^* a} = \{b: a_p = 0 \implies b_p = 0\}.
$$
\end{lemma}

\begin{proof} The inclusion $\subset$ is clear, and we only need to
show that if $b_p$ is zero whenever $a_p$ is zero, then $b$ can be
approximated  from within the orbit of $a$. Since $\adeles_f$ is the
restricted product of the $\Q_p$ relative to the $\Z_p$,  there is an
integer $n$ such that  $(n a)_p \in \Z_p$ for all $p\in \primes$, and
$n a$ has the  same orbit and the same zeros as $a$. So we may suppose
without loss of generality that
$a\in \prod_p \Z_p \times \R$.
The same argument shows that we might as well assume that
$b \in  \prod_p \Z_p \times \R$.

Suppose then that $b \in \prod_p \Z_p \times \R$ satisfies $b_p =0 $ if
$a_p = 0$. As a subspace of $\adeles$, the set $ \prod_p \Z_p \times \R$
has the product topology, so a typical basic open neighbourhood of $b$
has the form
$$ V = \prod_{p\in F} V_p \times  \prod_{p\notin F} \Z_p \times (x,y),
$$ where $F$ is a finite subset of primes and $(x,y)$ is an open
interval in $\R$.

Let $F' := \{p \in F: a_p \neq 0\}$; notice that whenever $p\in
F\setminus F'$ we trivially have
$r a_p = 0 = a_p = b_p \in V_p$ for all $r\in \Q^*$. We want to find
$r \in \Q^*$ such that 
\begin{itemize}
\item[\bf (i)] $r a_p \in V_p$ for $p \in F'$,
\smallskip
\item[\bf (ii)] $r a_p \in \Z_p$ for $p \in \primes \setminus F$,
\smallskip
\item[\bf (iii)]$r a_\infty \in (x,y)$.
\smallskip
\end{itemize}
For $p \in F'$, write $a_p$ as $p^{v_p(a_p)} u_p$ for some
unit $u_p \in \Z_p^*$, and choose  $k_p \in \Z$ such that $k_p \in
u_p\inv V_p$. Then $k_pp^{-v_p(a_p)} \in a_p\inv V_p$, and since
$a_p\inv V_p$ is open in $\Q_p$, there exists $\ell_p \in \N$ such that 
\begin{equation}\label{aVp} B\big(k_p{p^{-v_p(a_p)}} ,
p^{-\ell_p + v_p(a_p)} \big) \subset a_p\inv V_p.
\end{equation}
We will verify {\bf (i)} by finding $r$ in this ball.

 We deal with two cases
separately.

{\noindent \bf Case I: } Suppose there exists a prime $Q$ such that
$a_Q=0$.
If $Q=\infty$, the result follows from \lemref{binqa}. So we suppose
$Q\not=\infty$.

Let  $N_{F'} = \prod_{p\in F'} (p^{\ell_p - v_p(a_p)})$.  Choose $m \in
\N$ such that 
$$
\frac{N_{F'}}{Q^m} < \frac{y-x}{a_\infty}.
$$
By the Chinese Remainder Theorem, the congruences 
$$ n = k_p Q^m \left( \prod_{q\in F'\setminus \{p\}}
q^{v_q(a_q)}\right) \mod{p^{\ell_p}} 
\qquad  \text{ for } p\in F'
$$ have a solution $n_0$, and the set of solutions is then
$ n_0 + \left(\prod_{p\in F'} p^{\ell_p} \right) \Z$.
We claim that if $n$ is one of these solutions and we define $r$ by
\begin{equation}\label{defr} r := \frac{n}{Q^m \prod_{q\in
F'}q^{v_q(a_q)}},  
\end{equation} 
then $r$ is in the ball of \eqref{aVp}. To see this, we compute
\begin{eqnarray*}
\left| \frac{k_p}{p^{v_p(a_p)}} - \frac{n}{Q^m \prod_{q\in
F'}q^{v_q(a_q)}} \right|_p & =& 
\left| \frac{k_p Q^m \prod_{q\in F'\setminus\{p\}}q^{v_q(a_q)}- n }{Q^m
\prod_{q\in F'}q^{v_q(a_q)}}\right|_p \\ &=& \left| \frac{c
p^{\ell_p}}  {Q^m \prod_{q\in F'}q^{v_q(a_q)}}\right|_p 
\qquad \text{ for some $c\in\Z$}\\ &\leq& p^{-\ell_p + v_p(a_p)} 
\qquad \text{ by definition of $|\cdot|_p$ }.
\end{eqnarray*}
This says precisely that every $r$ of the form \eqref{defr} lies
in the ball of \eqref{aVp}, and hence satisfies {\bf (i)}.

To see that such $r$ also satisfy {\bf (ii)},  note that $a_Q =
0$ implies $r a_Q  = 0 \in \Z_Q$. If $p \notin F \cup \{Q\}$, then
$r \in \Z_p$ because  $p$ does not divide the denominator.
Hence $r a_p \in \Z_p$  for all $p\notin F$.

Now  recall that we chose $m$ large enough to ensure that every interval
of length $(y-x)/{a_\infty}$ contains at least one of the elements
of 
$$
\frac{ n_0 + \left(\prod_{p\in F'} p^{\ell_p} \right) \Z}{Q^m
\prod_{q\in F'}q^{v_q(a_q)}} =
\frac{ n_0 }{Q^m \prod_{q\in F'}q^{v_q(a_q)}} + \frac{N_{F'} \Z}{Q^m}.
$$ 
Therefore there exists $n \in n_0 + \left(\prod_{p\in F'} p^{\ell_p}
\right) \Z$ such that  the corresponding $r$ satisfies $r a_\infty \in
(x,y)$. This $r$ satisfies {\bf (i)}, {\bf (ii)} and {\bf (iii)}.

{\noindent \bf Case II: } Suppose $a_p \neq 0$ for every $p \in \primes
\cup \{\infty\}$. Then, because $a$ is not invertible, there must exist
infinitely many primes dividing
$a$.  Let 
$ N_{F} = \prod_{p\in F} (p^{\ell_p - v_p(a_p)})
$.
We can choose a
set $G$ of primes disjoint from $F$ such that every
$q\in G$ divides $a$, and such that
$$
\frac{N_{F}}{\prod_{q \in G} q} < \frac{y-x}{a_\infty}.
$$ 
A similar argument to that of Case I, with $Q^m$ replaced  by
$\prod_{q\in G} q$, now yields a rational number $r$
satisfying {\bf (i)}, {\bf (ii)} and  {\bf (iii)}.
\end{proof}

\begin{proposition}\label{bijection}
The quasi-orbit map $q : \adeles \to \mathcal Q(\adeles/\Q^*)$ is given by:
$$
a\in \adeles \mapsto q(a) = 
						\begin{cases} 
         \Q^*a & \text{ if $a$ is invertible},\\
         \{b\in \adeles \setminus \adeles^*:  b_p = 0 \text{ iff } a_p =0\} &
           \text{ if $a$ is not invertible}.	\end{cases}
$$
\end{proposition}
\begin{proof}
The result is a direct application of \lemref{qorb-inv} and \lemref{qorb-noninv}.
\end{proof}

We aim to obtain a better description of this quasi-orbit space and of its topology,
along the lines of \proref{qospace}. This will involve parametrising the
quasi-orbit space and then describing the topology on the parameter space.

One part of the parameter space will be the power set $2^\PP$ of the set of
extended primes  $\PP := \primes \cup \infty$;  we again
endow $2^\PP$ with its power-cofinite topology, generated by the basic open sets
$U_G = 2^{\PP\setminus G}$, where $G$ runs through the finite subsets of $\PP$. 

The other part of the parameter space will be $\mathcal U :=\prod_{p \in \primes}
\Z_p^* \times \R^*_+$ with the product topology, which is the topology of $\mathcal U$ as a
subspace of the ideles, see for instance \cite[Chapter IV.3]{weil} or \cite[\S 5.1]{ram-val}.
Although the topology of the ideles
as a restricted product is strictly stronger than the subspace topology 
they inherit from $\adeles$, the two topologies coincide on $\mathcal U$.
 To see this, suppose $u_n \to u$ in the adelic topology of $\mathcal U$. Then for each $p$,
$(u_n)_p$ converges to $(u)_p$ in $\Z_p$, and hence in $\Z_p^*$, because
$u$ and $u_n$ all lie in the subset $\mathcal U$ of $\adeles$. It follows that $u_n \to u$ in
the product topology of $\mathcal U$, that is, in the idelic topology.

\begin{proposition}
The map $ \chi: \adeles  \to  2^\PP \sqcup \mathcal U$ defined by
$$
a \in \adeles \mapsto \chi(a) = 
\begin{cases}
            a \sign(a_\infty) \prod_{\primes} p^{-v_p(a)} \in \mathcal U & \text{ if $a$ is
            invertible},\\
            \{p \in \PP: a_p =0\} \in 2^\PP & \text{ if $a$ is not invertible}.
\end{cases}
$$ 
factors through the quasi-orbit map, and induces a 
bijection 
$ \chi(a) \mapsto  q(a)$ of $2^\PP \sqcup\mathcal U$ onto the quasi-orbit space $\mathcal
Q(\adeles/\Q^*)$.
\end{proposition}
\begin{proof}
Every subset of $\PP$ is the vanishing set of a noninvertible adele,
and every $u\in \mathcal U$ is equal to $\chi(u)$, so $\chi$ maps $\adeles$ onto $2^\PP \sqcup\mathcal U$.
It remains to show that $\chi(a) = \chi(b)$ if and only if 
$q(a) = q(b)$ for every $a$, $b$ in $\adeles$.

 Suppose $\chi(a) = \chi(b)$. This implies that either $a$ and $b$ are both invertible,
or else they are both not invertible.
\begin{itemize}
\item If they are both invertible, we have 
$$
a \sign(a_\infty) \prod_{\primes} p^{-v_p(a)} =
b \sign(b_\infty) \prod_{\primes} p^{-v_p(b)},
$$ hence  $b = \pm a \prod_{\primes} p^{-v_p(a) +
v_p(b)} \in \Q^* a$, from which  $q(a) = q(b)$.

\item If they are both not invertible, then 
$\{ p \in \PP: a_p=0 \} = \{ p\in \PP:b_p =0 \}$ and
clearly $q(a) = q(b)$.
\end{itemize}
Suppose now $q(a) = q(b)$. This also implies that either $a$ and $b$ are both invertible,
or else they are both not invertible.

\begin{itemize}
\item If they are both invertible, we can write $b = r a$ for some $r\in \Q^*$, and then
\begin{align*}
\chi(b) &= (ra) \sign(ra_\infty)\prod_{\primes} p^{-v_p(ra)} \\
&= (r \sign(r)\prod_{\primes} p^{-v_p(r)}) (a \sign(a_\infty)\prod_{\primes} p^{-v_p(a)}) \\
&= a \sign(a_\infty)\prod_{\primes} p^{-v_p(a)} 
= \chi(a),
\end{align*}
 because
$r \sign(r) \prod_{\primes} p^{-v_p(r)} = 1$. 

\item If they are both not invertible, then 
\[\{c\in \adeles \setminus \adeles^*:  c_p = 0 \text{ iff } a_p =0\} =
\{c\in \adeles \setminus \adeles^*:  c_p = 0 \text{ iff } b_p =0\},\]
and from this it follows that $a_p = 0$ if and only if $b_p =0$, so $\chi(a) = \chi(b)$.
\end{itemize}
This completes the proof.
\end{proof}

The set $\mathcal U$ is a locally compact group under multiplication, and the restriction of 
$\chi$ to the ideles $\adeles^*$ is a group homomorphism onto $\mathcal U$ with kernel
$\Q^*$. Thus $\mathcal U$ is isomorphic to the idele class group $\adeles^*/
\Q^*$.

The next step is to topologise $2^\PP \sqcup \mathcal U$ so as to make the bijection 
from \proref{bijection} a
homeomorphism. The right topology 
comes from the power-cofinite topology on $2^\PP$ and the product topology on $\mathcal U$,
but it is necessary to specify how these two parts interact. 
First, to deal with a subtlety arising from the difference between adelic and idelic
topologies, we need to consider the absolute value function on $\adeles$. 

The absolute value $\| a\|$ of an adele
$a$ is, by definition, the product of the normalized $p$-adic absolute values:
$$
\|a\| := |a_\infty| \prod_{p\in \primes} p^{-v_p(a)},
$$ 
where $v_p$ is the $p$-adic valuation and $|a_\infty|$ is the
usual absolute value of the real component of $a$. The absolute
value vanishes at $a$ if and only if either $a_p =0$ for
some $p$ or $v_p(a) > 0$ for infinitely many $p$; in other words, 
$\|a\| = 0 $ if and only if  $a$ is not
invertible. Since noninvertibles are dense in $\adeles$, the absolute value is not a continuous
function.

\begin{lemma}\label{uppersc}
The absolute value is an upper semi-continuous function on $\adeles$.
\end{lemma}
\begin{proof}
The maps $a \mapsto p^{-v_p(a)}$ for $p\in \primes$ and $a\mapsto |a_\infty|$ are all continuous
on $\adeles$, so for each finite subset $F$ of $\primes$ the product  
$$
\xi_F: a \in \adeles \mapsto |a_\infty| \prod_{p \in F} p^{-v_p(a)}
$$
is continuous, and, the net $(\xi_F )_F$ directed by the finite subsets of $\primes$ under inclusion
converges pointwise to the absolute value.
For each finite $G \subset \primes$ and each choice of $k_p \in \Z$ for $p\in G$, the set
$$
W := \prod_{p\in G} p^{k_p} \Z_p \times \prod_{p\notin G} \Z_p \times \R
$$
is open in $\adeles$, and every adele belongs to some $W$
(given $a\in \adeles$, it suffices to choose $G = \{p: v_p(a) < 0\}$ and $k_p = v_p(a)$).
Moreover, for every $b\in W$, if $v_p(b) < 0$ then $p \in G$. It follows that
the tail net $\{\xi_F\}_{F \supset G}$ is nonincreasing when restricted to $W$, because
all the factors greater than $1$ have been included in $\xi_F(b)$ already.
Thus, the restriction of $\|\cdot \| = \lim_F \xi_F$ to $W$ is the infimum of the restricted tail net
and hence is upper semicontinous on $W$. Since each $W$ is open and their union is all of $\adeles$,
this completes the proof.
\end{proof}

From \eqref{factor-ideles} we know that the map $(r,u) \mapsto ru$ is a bijection of $\Q^*
\times \mathcal U$ onto $\adeles^*$. Since $\Q^*$ has the discrete topology, and the product (idelic)
topology on $\mathcal U$ coincides with the adelic topology, it is clear that the map is continuous.
However, its inverse is not continuous: let $p_n$ be the $n$th prime and define a
sequence $a_n \in \adeles^*$ by $(a_n)_q = 1$ for $q \neq p_n$ and
$(a_n)_{p_n} = p_n$; then $a_n \to 1\in \adeles$, but $(u_n)_\infty \to 0$, so
$(r_n,u_n)$ cannot converge in $ \Q^* \times \mathcal U$. Notice that the absolute value $\|a_n\|$
tends to zero in this example; the following lemma shows that this is crucial.

\begin{lemma}\label{F=UXQ}
For each $\varepsilon > 0$ let $F_\varepsilon := \{a\in \adeles: \|a\| \geq \varepsilon\}$
and $\mathcal U_\varepsilon := \{ u\in \mathcal U: u_\infty \geq \varepsilon\}$,
considered with their relative topologies as subsets of the adeles and the ideles, respectively.
Then $F_\varepsilon$ is closed, $\chi(F_\varepsilon) = \mathcal U_\varepsilon $, and the map $(r,u)
\mapsto ru$ is a homeomorphism of $ \Q^* \times \mathcal U_\varepsilon  $ onto $F_\varepsilon$. 
\end{lemma}

\begin{proof}
By \lemref{uppersc}, $F_\varepsilon$ is closed, and 
the preceding discussion establishes that $(r,u) \mapsto ru$ is a continuous bijection, so
it only remains to show that the inverse is continuous. Suppose  $a_n\to a$ in $F_\varepsilon$; we
need to show that
$u_i :=\chi(a_i)$  converges to $u:=\chi(a)$ in
$\mathcal U_\varepsilon $ and that $r_i:= r(a_i) $ is eventually equal to
$r(a)$.  Dividing everything by $r(a)$, we may assume that $a \in \mathcal U_\varepsilon$, and show
that $r_i$ is eventually $1$.  
For each $n \in \N$, let $p_n$ be the $n$th prime, and let
$$
U_n := \prod_{j \leq n} \Z_{p_j}^* \times \prod_{j > n} \Z_{p_j} \times I
$$
where $I$ is an interval in $\R^*_+$ containing $a_\infty$. 
For each $n$,  $a_i $ is eventually in $U_n$, and then
$$
v_{p_j}(r_i) = 0 \ \text{ for } j \leq n \quad \text{ and } \quad
v_{p_j}(r_i) \geq 0 \  \text{ for } j > n.
$$
It follows that $r_i$ is eventually an integer, and that 
either $r_i = 1$ or $r_i \geq p_i$ 
(because $p_i$ is the smallest possible prime factor of $r(a_i)$).
Since $ r_i u_i = a_i$ converges to $a \in F_\varepsilon$,
we know that $(r_i u_i)_\infty $ converges to $a_\infty \geq \varepsilon$.
But $(u_\lambda)_\infty$ is bounded away from $0$, so $\{r_i\}$ has to remain bounded,
and this can only happen if $r_i$ is eventually $1$. 
\end{proof}

\begin{theorem}
Let $\tau$ be the topology on the parameter space $2^\PP \sqcup \mathcal U$ pulled back from the
quotient topology on the quasi-orbit space under the bijection of \proref{bijection}
(so that $\tau$ is the quotient topology induced by the map
$\chi:\adeles \to 2^\PP \sqcup \mathcal U$).
Then the $\tau$-closure of a set $B \subset 2^\PP$ is 
$$
\overline{B}^\tau = \begin{cases} \text{the power-cofinite closure of } B & \text{ if $B$ is not
power-cofinite dense in } 2^\PP,\\
2^\PP \sqcup \mathcal U & \text{ if $B$ is  power-cofinite dense in } 2^\PP,
\end{cases}
$$
and the $\tau$-closure of a set $C \subset \mathcal U$ is
$$
\overline{C}^\tau = \begin{cases} \text{the idelic closure of $C$ } & \text{ if $\{c_\infty: c\in
C\}$ is bounded away from $0$,}\\
 2^\PP \sqcup \mathcal U & \text{ if $0 \in \overline{\{c_\infty: c\in C\}}$ .}
\end{cases}
$$
\end{theorem}
\begin{proof}
First we show that the power-cofinite closure of $B$ is always contained in 
$\overline{B}^\tau$. Let $T$ be in the power-cofinite closure of $B$;
we will show that that every (basic) open set
containing $T$ intersects $B$. We can assume that the basic open set is  
$\chi(W)$ for $W = \prod_F V_p \times \prod_{F^c} \Z_p$, where $V_p$ is an open subset of $\Q_p$ for
every $p$ in the finite subset $F$ of $\PP$, and that $T =\chi(t)$ for $t\in W$.
We need  to find an adele $w$ such that
\begin{itemize}
\item  $w_p \in V_p$ for $p\in F$ and $w_p\in \Z_p$ for $p\notin F$ (so that $w\in W$);
\item $w$ is not invertible in $\adeles$ (so that $\chi(w) = \{p: w_p=0\}$);
\item  $\chi(w) \in B$.
\end{itemize}
The set $G = \{p\in F: 0 \notin V_p\}$ is finite and disjoint from $T = \{p: t_p
=0\}$, because $t\in W$. Thus $T\in U_G$, and because $T$ is in the power-cofinite closure, there
exists
$b\in B
\cap U_G$ --- that is,
$b\cap G
=\emptyset$. Choose $w$ as follows:
$$
w_p = \begin{cases} t_p & \text{ if $p\in G$ }\\
                        p & \text{ if $p \notin b\cup G$}\\
                        0 & \text{ if $p \in b$}.
  \end{cases}
$$ 
Then $w$ is not invertible (even if $b=\emptyset$) and $\chi(w) = b$. This proves that $T \in
\overline{B}^\tau$.

Next we show that if $B$ is not power-cofinite dense in $2^\PP$ then 
 the power-cofinite closure $\overline{B}^{pc}$ of $B$ contains $\overline{B}^\tau$.
Because $B$ is not dense, it misses some basic open set $U_G$, and $\overline{B}^{pc}\subset
2^\PP\setminus U_G$. We claim that $2^\PP\setminus U_G$ is also $\tau$-closed. 
To see this,  write 
$$
2^\PP \setminus U_G = \chi(\chi\inv(2^\PP \setminus U_G)) 
= \chi(\{a\in \adeles: \chi(a)\cap G \neq \emptyset\}).
$$
and then observe that the set
\begin{eqnarray*}
 \{a\in \adeles: \chi(a)\cap G \neq \emptyset\} 
   &=&  \{a\in \adeles: \{p: a_p =0\} \cap G \neq \emptyset\} \\
   &=&  \{a\in \adeles:  a_p =0  \text{ for some } p\in G \} \\
   &=&  \bigcup_{p\in G} \{a \in \adeles: a_p =0\}
\end{eqnarray*}
is  a finite union of $\Q^*$-invariant closed
sets. Thus $2^\PP\setminus U_G=\chi(\{a\in \adeles: \chi(a)\cap G \neq \emptyset\}) $ is
$\tau$-closed, as claimed.
This implies that the $\tau$-closure of $B$ is contained in $2^\PP$, and indeed is
contained in
\[
\overline{B}^{pc}=\cap\{2^\PP\setminus U_G: B\cap U_G=\emptyset\}.
\]

If $B$ is power-cofinite dense in $2^\PP$, then $\emptyset$ is in the power-cofinite closure of
$B$.  The initial paragraph of the proof, with $T = \emptyset$, shows that $\emptyset \in
\overline{B}^\tau$; we will show that $\{\emptyset\}$ itself is $\tau$-dense.  Let $W$ be any open
set in $\adeles$. Choose a non-invertible adele $a$ such that $a_p\not=0$ for every $p$. 
\lemref{qorb-noninv} implies that
$W$ contains $ra$ for some $r\in\Q^*$. 
Then $\emptyset = \chi(a) =\chi(ra) \in \chi (W)$. Thus $\{\emptyset\}$ is $\tau$-dense, and we must
have
$\overline{B}^\tau =  2^\PP \sqcup \mathcal U$. This completes the description of
$\overline{B}^\tau$.

Next let $C\subset \mathcal U$, and assume first that $0 \in \overline{\{c_\infty: c\in C\}}$. We
claim that every singleton $\{p\}$ belongs to
$\overline{C}^\tau$. There is a strictly increasing sequence $n_k$ such that each interval
$(p^{-(n_k + 1)}, p^{-n_k}]$ contains $c_\infty$ for some $c\in C$, and for each $k$ we choose one
such $c = c_k$.  Since $n_k \to \infty$, we have $(p^{n_k+1} c_k )_p \to 0$ in $\Z_p$. Moreover, for
all $k$ we have $(p^{n_k+1} c_k )_q \in \Z_q^*$ for every finite prime $q\neq p$ and
$(p^{n_k+1} c_k )_\infty \in (1,p]$. 
Since $\Z_p \times \prod_{q\neq p} \Z_q^*\times [1,p]$ is compact, we can assume by passing to a
subsequence that $(p^{n_k+1} c_k )_q \to a_q$ for every $q\in \PP$; then
$a_q$ vanishes precisely when $q=p$, and hence $\chi(a) = \{p\}$. By continuity of $\chi$, we
know that
$\chi(p^{n_k+1} c_k)\to\chi(a) =\{p\}$ in $(2^\PP \sqcup\mathcal U,\tau)$, and, since
$\chi(p^{n_k+1} c_k) = \chi(c_k) = c_k \in C$, it follows that $\{p\} \in
\overline{C}^\tau$, as claimed

Since every singleton $\{p\}$ is in $\overline{C}^\tau$, $\overline{C}^\tau$ meets every basic open
neighbourhood $U_G$ of $\emptyset$ in $2^\PP$, and $\emptyset\in\overline{C}^\tau$. Since we have
already seen that $\{\emptyset\}$ is $\tau$-dense, it follows that $\overline{C}^\tau=2^\PP\sqcup
\mathcal{U}$. 

Assume now  that
$c_\infty \geq \varepsilon$ for every $ c\in C$, and let 
$F_\varepsilon := \{a\in \adeles: \|a\| \geq \varepsilon\}$. 
The inclusion map $\mathcal U_\varepsilon \hookrightarrow F_\varepsilon \stackrel{\chi}{\to} 2^\PP
\sqcup \mathcal U$ is continuous, so the closure $\overline{C}^{\mathcal U}$ of $C$ in $\mathcal U$
is certainly contained in the $\tau$-closure of $C$ in $2^\PP \sqcup \mathcal U$.
On the other hand, because the inverse of $(r,u) \mapsto ru$ is continuous by \lemref{F=UXQ},
the set $\Q^* \overline{C}^{\mathcal U}$ is closed in $F_\varepsilon$. Since it is also 
$\Q^*$-invariant, its image $\overline{C}^{\mathcal U} = \chi(\Q^* \overline{C}^{\mathcal U})$
is closed in the quotient topology; in other words, it is $\tau$-closed. Thus 
$\overline{C}^\tau \subset \overline{C}^{\mathcal U}$, and we conclude that 
$\overline{C}^\tau = \overline{C}^{\mathcal U}$.
\end{proof}

From here it is relatively easy to describe $\Prim(C_0(\adeles) \rtimes \Q^*)$
and its topology.
Recall that for each $a\in \adeles$ there is a covariant representation $(M_a,
\lambda)$ of $(C_0(\adeles),\Q^*)$ on $\ell^2(\Q^*)$ given by:
$$ (M_{a} f )\xi (r) = f(ra) \xi(r) 
\quad \text{ and }\quad 
\lambda_s \xi (r)= \xi (s\inv r).
$$ 
The primitive ideals of $C_0( \adeles) \rtimes \Q^*$ 
are the kernels of these representations; there are three types
corresponding to the three types of orbit closures, and our description of the topology on 
$\Prim(C_0(\adeles) \rtimes \Q^*)$ follows from Williams' theorem (see the comments in 
Remark~\ref{disjunion}). 

\begin{corollary}\label{Primfulladeles}\ 

\begin{enumerate}
\item[(1)] The map $\gamma \mapsto \ker (\varepsilon_{0} \times \gamma)$ is
an injection of $\widehat{\Q^*}$ into $\Prim (C_0( \adeles) \rtimes
\Q^*)$, and
 its image is a closed subset of $\Prim (C_0( \adeles) \rtimes \Q^*)$.

\item[(2)] If $a \neq 0$ is not invertible, then $\ker(M_a\times \lambda)$
depends only on $S(a) := \{p \in \PP: a_p = 0 \}$, and the map 
$$ 
S(a) \mapsto \ker (M_a\times \lambda)
$$ 
is an injection of $ 2^\PP \setminus \{\PP\}$ into  $\Prim (C_0( \adeles) \rtimes \Q^*)$. 

\item[(3)]  If $a \in A$ is invertible, then $\ker(M_a\times \lambda) $
depends only on the orbit $\Q^* a$, which contains a unique element 
$u\in\mathcal U := \prod_{p \in \primes} \Z_p^* \times \R^*_+$; the map
$$ 
u \mapsto \ker(M_u\times \lambda)
$$ 
is an injection of the set $\mathcal{U}$ into  $\Prim (C_0( \adeles) \rtimes \Q^*)$. 
\end{enumerate} 
Taken together, these maps parametrize the primitive spectrum of 
$C_0( \adeles) \rtimes \Q^*$ in the sense that their images are
disjoint and their union is all of $\Prim (C_0(\adeles) \rtimes \Q^*)$.
Define $\overline{A}$ for a nonempty subset $A$ of $\widehat{\Q^*} \sqcup (2^\PP \setminus
\{\PP\}) \sqcup \mathcal U$  by
\begin{equation*}
\overline{A} := \begin{cases} 
\overline{A}^{\widehat{\Q^*}} & \text{ if }A \subset \widehat{\Q^*};\\
\widehat{\Q^*}\sqcup \big(\overline{A}^{pc} \setminus \{\PP\} \big)  & 
\text{ if }  A \subset 2^\PP \setminus \{\PP\} \text{ and } \overline{A}^{pc} \neq
2^\PP;\\
\overline{A}^{\mathcal U} & 
\text{ if } A \subset \mathcal U \text{ and } 0 \notin \overline{\{\|a\|: a\in A\}};\\
\widehat{\Q^*} \sqcup (2^\PP \setminus \{\PP\}) \sqcup \mathcal U &
\begin{cases} \text{ if } A\subset  2^\PP \setminus \{\PP\} \text{ and } \overline{A}^{pc} = 2^\PP &
\\
    \text{ or if }
 A \subset \mathcal U \text{ and } 0 \in \overline{\{\|a\|: a\in A\}}. &  
\end{cases} \end{cases}
\end{equation*} 
Then $A \mapsto \overline{A}$ (has an obvious extension that) satisfies Kuratowski's closure axioms,
and the resulting topology on $\widehat{\Q^*} \sqcup (2^\PP \setminus \{\PP\}) \sqcup \mathcal U$ 
makes the parametrisation a homeomorphism.
\end{corollary}

\begin{remark}
(1)  The injection of $\widehat{\Q^*}$ is a homeomorphism onto its image in
$\Prim(C_0(\adeles)\rtimes\Q^*)$.  Indeed, the map
$\varepsilon_0$ induces an isomorphism of 
$(C_0( \adeles) \rtimes \Q^*) / ((\ker \varepsilon_0) \times \Q^*)$ onto
$C^*(\Q^*)$, and $\widehat{\Q^*}$ is the primitive ideal space of the
quotient viewed as a closed subset of $\Prim(C_0( \adeles) \rtimes
\Q^*)$.

(2) The second injection is also a homeomorphism onto its image: to see
this, just note that the power-cofinite closure of $A\subset 2^\PP$ is
equal to $\overline{A}^\tau \cap 2^\PP$.

(3) The third injection is not a homeomorphism onto its range.  Indeed, the image under
$\chi$ of the sequence
$\{a_n\}$ defined before \lemref{F=UXQ} is closed in the idelic topology on $\mathcal{U}$, but is
$\tau$-dense in $2^\PP\sqcup \mathcal{U}$, and hence dense for the relative topology on
$\mathcal{U}\subset 2^\PP\sqcup\mathcal{U}$.

(4) The representations associated to elements of $\mathcal U$ are
all CCR representations: the map $r \mapsto ru$ is a homeomorphism of
$\Q^*$ onto the discrete set $\Q^* u$, and hence the image of the representation 
$M_u\times \lambda$ is isomorphic to $\mathcal K (\ell^2(\Q^*))$ by the Stone-von Neumann theorem.
\end{remark}

\end{document}